\documentclass[onecolumn,aps,showkeys,nofootinbib]{revtex4}
\usepackage{epsfig}

\usepackage{natbib}
\usepackage{graphicx}

\usepackage{pdfpages}

\usepackage{lipsum}

\usepackage{amsfonts}

\usepackage{amsmath}
\usepackage{amssymb}
\usepackage{latexsym}
\usepackage{graphicx}
\usepackage{hyperref}
\usepackage{bm}

\usepackage{color,graphicx,subfigure,graphics,epsfig,amsmath,amssymb}
\usepackage{epsf}
\hypersetup{urlcolor=blue, colorlinks=true}  

\usepackage{lipsum}
\usepackage{hyperref}
\usepackage{bm}
\usepackage{enumerate}
\usepackage{amsthm}

\usepackage{bm}
\usepackage{enumerate}
\usepackage{amsthm}

\usepackage{stmaryrd}

\usepackage{systeme}

\def\p{{\cal P}}
\def\s{{\cal S}}

\def\t{{\cal T}}

\newtheorem{thm}{Theorem}[section]

\newtheorem{kor}[thm]{Corollary}

\theoremstyle{definition}
\newtheorem{defi}{Definition}

\theoremstyle{lemma}
\newtheorem{lemm}{Lemma}

\begin{document}
\input{epsf.sty}

\title{Novel Aspects Of The Global Regularity Of Primes}
\author{ Patrice M. Okouma$^{1}$ , Guillaume Hawing$^{2}$,  \\ 
\href{mailto:patrice.okouma@uct.ac.za}{patrice.okouma@uct.ac.za} \\
\it $^1$ Department of Mathematics \& Applied Mathematics, University of Cape Town, South Africa\\ 
\it  $^2$ Universit\'e Mahatma Gandhi, Conakry, R\'epublique de Guin\'ee \\
}
\begin{abstract}

A new deterministic scheme for characterizing the organization of primes is established. At its core are eleven generic tables, coupled with a three-criteria test applied on differences between pairs of consecutive odd composite numbers obtained from the tables. Our scheme trivially yields $ \p\left( m \right),$ the set of all primes below any preassigned limit $m$ in $ \mathbb{N}: \pi(m) = \left \vert  \p\left( m \right)  \right \vert.$ Setting $ m = 3 \left( 2n + 1 \right), \forall \ n \in \mathbb{N}^* ,$ we establish, for the first time, that $\pi(m) = 4 + \left \vert A_4(m) \right \vert + 2 \left \vert  A_6(m) \right \vert $ where  $ \left \vert  A_4(m)  \right \vert $ and $ \left \vert  A_6(m)  \right \vert$ are the cardinals of the sets of pairs of consecutive odd composite numbers, below $m$ and having, each, a difference of $4$ and $6,$ respectively. {\it We prove that $ A_6(m) $ is the set of Twin Primes below $m.$} $A_4(m) $ is the set of primes below $m$ that we call "isolated" as they each appear, alone, between a pair of consecutive odd composite numbers distant by $4$. $\t(m)$ is the number of consecutive odd composite numbers below $m$. We establish that $\t(m) = 1 + \left \vert A_2(m) \right \vert + \left \vert A_4(m) \right \vert + \left \vert  A_6(m) \right \vert$ where $ \left \vert  A_2(m)  \right \vert $ is the number of pairs of consecutive odd composite numbers, below $m$, distant by $2$: $ \left \vert  A_2(m)  \right \vert  = \frac{m - 9}{2} - 3 \left \vert  A_6(m)  \right \vert - 2 \left \vert  A_4(m)  \right \vert .$ With {\it $m$ an odd number}, combining $\pi(m)$, $ \t(m),$ and $ \left \vert  A_2(m)  \right \vert, $ we establish that  $$ 4 \left \vert  A_6(m)  \right \vert + 7  = m -  2 \left(  \t(m) +  \left \vert  A_4(m)  \right \vert \right)  > 0, \ hence, \ 0 < 1 -  \frac{ 2 }{m} \left(  \t(m) +  \left \vert  A_4(m)  \right \vert \right)  < 1 , \ \mbox{for any $m$ as above.} $$ The non empty and unique set $\s \equiv \left\lbrace 1 -  \frac{ 2 }{l} \left(  \t(l) +  \left \vert  A_4(l)  \right \vert \right) \ \mbox{such that} \ \mbox{l = 3(2k + 1), $k \in \mathbb{N}^*$} \right\rbrace \subset \mathbb{R} $ is bounded and independent of m. 
Hence $$ Inf \left( \s \right)  \ \mbox{{\it exists, is unique, finite and independent of m.}}$$ 
By definition, $ Inf \left( \s \right) > 0.$ $\mathbb{Q}$ Dense in $\mathbb{R}$ also guarantees the latter. We then introduce
\begin{equation} 
\alpha = Inf \left( \s \right) \ \mbox{independent of m.} 
\end{equation}
We have $$ 0 < \alpha \leq 1 -  \frac{ 2 }{m} \left(  \t(m) +  \left \vert  A_4(m)  \right \vert \right) < 1 , \ \mbox{for any m = 3(2n + 1), $ n \geq 1$, an integer.}$$
In others words, $$ m - 2 \left( \t(m) +  \left \vert  A_4(m)  \right \vert \right) = 4 \left \vert  A_6(m)  \right \vert + 7 \geq \alpha m \Leftrightarrow 2 \left \vert  A_6(m)  \right \vert \geq \frac{\alpha m - 7}{2} , \ \mbox{ with $ 1 > \alpha > 0$ and m as above.} $$
{\it The last relation above establishes, for the first time, that} $$ \displaystyle \lim_{ m \longrightarrow \infty }  2 \left \vert  A_6(m)  \right \vert  = + \infty. $$ 
Similarly we establish that, $$ \displaystyle \lim_{ m \longrightarrow \infty }  \left \vert  A_4(m)  \right \vert  = + \infty.$$
Thanks to $\pi(m),$ we also have a new proof of Euclid's Theorem. Through $\t(m),$ there is an infinite number of odd composite numbers.

\end{abstract}

 
\maketitle

\section{Introduction}

Prime numbers play an important role, both in pure mathematics and its applications. Primes occur in a very irregular way within the sequences of natural numbers. In particular, the distribution of prime numbers exhibits a local irregularity but a global regularity. One of the best existing results illustrating that global regularity is the prime number theorem giving the  number of prime number, $\pi(n)$, not exceeding any upper limit $ n \in  \mathbb{N} $. In an epoch-marking work published in the late nineteenth century, it was shown that $\pi(n)$ behaves asymptotically as $\frac{n}{\ln n}$ \cite{Poussin}. Prime numbers derive most of their importance from the fundamental theorem of arithmetic \cite{Gauss1798}. Given their importance in the construction of all natural numbers, a problem that presents itself at the very roots of mathematics is therefore the one about the organization of the primes among the integers. In this work, using consecutive odd composite numbers as a gateway, we present a new deterministic scheme allowing an unexpected foray into the organization of prime numbers, including a novel re-investigation of the Twin Primes conjecture. We recall that Twin Primes are pairs of Prime numbers separated by two. Some attribute the Conjecture to the mathematician Euclid of Alexandria. The first formulation of the Twin Primes conjecture, in its modern form, was given in 1846 by French mathematician Alphonse De Polignac \cite{dePolignac1849}. The conjecture states that there are infinitely many Twin Primes. In April 2013, Zhang published his results, showing that there were infinitely many primes with a gap less than 70 million \cite{Zhang2014}. Working independently of Zhang and the Polymath Project, Maynard has established a bound of 600 using an entirely different method \cite{Maynard2019}. In April 2014, within the Polymath Project, the bound stood at 246 \cite{Polymath2014a}. Assuming an additional result on the distribution of primes, known as the Elliott-Halberstam conjecture, the bound can be reduced to 6 \cite{Polymath2014b}. The present paper unveils some novel aspects of the organization of Primes. It is a first in a series. In section {\bf II}, we comprehensively present the foundations of our new Deterministic scheme. In  section {\bf III}, we present a few consequences of the scheme, including a new light on the Twin Primes Conjecture. We conclude in section  {\bf IV}.  

\section{A Novel Deterministic Scheme Yielding the Set of Primes up to any $ n \in  \mathbb{N} $}

Towards clarifying our approach and some terminology, let us start with some definitions. 
\begin{defi}
Any $ \ p > 1, \ \mbox{in} \ \mathbb{N},$ is an odd composite number if and only if
\begin{enumerate}[(1)]
 \item $p$ is an odd number;
 \item $p$ is not a prime number. 
 \label{def_con}
\end{enumerate}
\end{defi}
We note that $1$ is neither prime or composite. 
\begin{defi}
For any $ n \ and \ p \ elements \ of \ \mathbb{N}, \ $ $n$ and $p$ are consecutive odd composite numbers (COCOONs) if and only if there exists no other odd composite number $k$ such that $ n < k < p.$
\label{def_cocon}
\end{defi}

By construction, any odd (composite) number has to end with $1 \ or \ 3 \ or \ 5 \ or \ 7 \ or \ 9.$ Hence,
\begin{lemm}
\begin{enumerate}[] For some $n \in \mathbb{N^*}$
\item $PN_1 \equiv 10n + 1 \ ,$ for all odd numbers ending with $1$;
\item $PN_3 \equiv 10n + 3 \ ,$ for all odd numbers ending with $3$;
\item $PN_5 \equiv 10n + 5 \ ,$ for all odd numbers ending with $5$;
\item $PN_7 \equiv 10n + 7 \ ,$ for all odd numbers ending with $7$;
\item $PN_9 \equiv 10n + 9 \ ,$ for all odd numbers ending with $9$.
\end{enumerate}
\label{struct}
\end{lemm}

The fundamental theorem of arithmetic establishes that any odd composite number can only be written as a product of odd numbers: 
\begin{lemm} 
 For any $ \mbox{odd composite number, $ p$ , there exists} \ l \  \mbox{and}  \ m \mbox{ each an odd number such that } p = l*m .$ 
\label{on_decomp}
\end{lemm}

Based on lemma \ref{struct}, one can express all possible products $p=l*m$ for $p$ any odd number. This consists in writing down all possible products $ PN_i * PN_j $ for $i,j  \in \left\lbrace 1, 3, 5, 7, 9  \right\rbrace,$ the only options for the last digit of each $PN_k$. All potential products yield fifteen (15) equalities. Five redundant equalities associated to all odd composite numbers ending with 5 are collapsed into a single relation. This leaves us with a set of eleven (11) equations defining all possible odd composite numbers. The eleven equations consist in three (for any $p$ ending with 1), two (3), one (5), two (7) and three (9). These eleven relations can be formulated into eleven (11) tables. For illustration purposes, the two tables yielding all odd composite numbers ending with 7 are given in Tables \ref{con7_a} and \ref{con7_b}. From these eleven relations, when computing the complete list of odd composite number below any integer $n \geq 9,$  the only constrain is that each product $p$ from each of the eleven tables has to be smaller or equal to $n.$ Once ordered, one gets the set of consecutive odd composite numbers below $n.$ All odd numbers not in that final set have to be non composite i.e. primes. Adding $\left \lbrace 2, 3 , 5 , 7 \right \rbrace$ completes  $\mathcal{P}(n),$ the set of all primes of below $n \geq 9$: $\pi(n) \equiv \left \vert \mathcal{P}(n) \right \vert.$ Theorem \ref{decision_test}, our first main result, underlies our new Deterministic scheme. 

\begin{table}[h]
\centering
\caption{Table representation for $ \left\lbrace  \left( 10 k_1 + 1 \right) \left(10 k_2 + 7 \right ), \ k_1 \in \mathbb{N}^* \ and \ k_2 \ \in \mathbb{N} \right \rbrace $ : odd composite numbers ending with 7}
\begin{tabular}{|c|c|c|c|}\hline
$k_1 $ & $ k_2 $ & \ $ \left( 10 k_1 + 1 \right)   $ & $   \left(10 k_2 + 7 \right) $ \\ \hline
\ 1 & 0 & 11 & 7 \\ \hline
\ 2 & 1 & 21 & 17 \\ \hline
\ ..& ... & ... & ... \\ \hline
\ $\infty$ & $\infty$ & $\infty$ & $\infty$ \\ \hline
\end{tabular}
\label{con7_a}
\end{table}

\begin{table}[h]
\centering
\caption{Table representation for $ \left \lbrace \left( 10 k_3 + 3  \right) \left( 10 k_4 + 9\right), \ k_3, \ k_4 \ \in \mathbb{N}   \right \rbrace $:  odd composite numbers ending with 7. Together with Table \ref{con7_a} they yield all odd composite numbers ending with 7. }
\begin{tabular}{|c|c|c|}\hline
$k_3 , k_4 $ & \ $ \left( 10 k_3 + 3  \right)   $ & $  \left( 10 k_4 + 9\right) $ \\ \hline
\ 0 & 3 & 9 \\ \hline
\ 1 & 13 & 19 \\ \hline
\ 2 & 23 & 29 \\ \hline
\ ... & ... & ... \\ \hline
\ $\infty$ & $\infty$ & $\infty$ \\ \hline
\end{tabular}
\label{con7_b}
\end{table}

\begin{thm}[Main Theorem]

 For any pair $ \left \lbrace N_1 \ ,  N_2 \right \rbrace $ of consecutive\footnote{consecutivity as in definition 2 is crucial.} odd composite numbers,
 \begin{enumerate}[(1)]
  \item $N_2 - N_1 = 2$ or $N_2 - N_1 = 4$ or $N_2 - N_1 = 6,$
  \item If $N_2 - N_1 = 2$ then there is no prime number between $N_1$ and $N_2$,
  \item If $N_2 - N_1 = 4$ then there is one prime numbers between $N_1$ and $N_2$,
  \item If $N_2 - N_1 = 6$ then there are two prime numbers between $N_1$ and $N_2$.
 \end{enumerate}
\label{decision_test}
\end{thm}

\begin{figure}[tbp]
 \includegraphics[width=1.\linewidth]{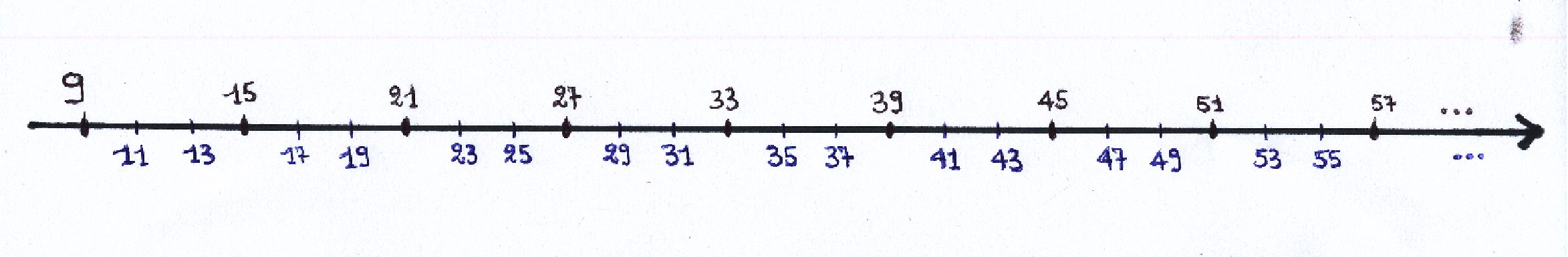}
 \caption{Partially illustrated here is an ordered set of all consecutive odd composite numbers (cocoons) each multiple of $3$ (in black) that are smaller or equal to some $n \in \ \mathbb{N}$. By construction, $9$ is the first odd composite number. $ \varOmega_{cocoon} \left( n \right) $ is the set of all cocoons  smaller or equal to $n.$  $ \varOmega_{cocoon} \left( n \right) \subseteq  \displaystyle \bigcup^m_{k = 1, \  k \in \mathbb{N} } \varOmega_k  $  with  $ \varOmega_k = \llbracket 3 \left( 1 + 2 k   \right)  \ .. \   3 \left( 3 + 2 k   \right) \rrbracket , k \in  \mathbb{N}^*  $  $ \subseteq \mathbb{N},$ ordered. The integers in Black are boundaries of {\it some} $ \varOmega_k $.}
 \label{fig:line_of_numbers}
 \end{figure}

\begin{proof} [Proof of Theorem \ref{decision_test}]

 \begin{enumerate} [(A)]
 \item We prove item (1) in Theorem \ref{decision_test}. Let us introduce the ordered set $\varOmega \left( n \right) $ of all odd composite numbers smaller or equal to $n \in \ \mathbb{N}.$ By construction, $9$ is the first element of $\varOmega \left( n \right).$  $ \varOmega \left( n \right) \subseteq  \displaystyle \bigcup^m_{k = 1, \  k \in \mathbb{N} } \varOmega_k  $  with $ \varOmega_k =  \left \lbrace \mbox{$p$ odd composite number such that $3 \left( 2 k + 1  \right) \leq p  \leq   3 \left( 2 k + 3  \right) , $ $ k \in \mathbb{N}^*  $ }  \right \rbrace $ $ \subseteq \mathbb{N}$. %
 $ \varOmega \left( n \right) $ can be visualized as a line of integers partially illustrated in Figure \ref{fig:line_of_numbers}.  
 \begin{itemize}
 \item We are interested in the separation between the elements of any pair of consecutive odd composite numbers within $ \varOmega \left( n \right) $; 
 \item As illustrated in Figure \ref{fig:line_of_numbers}  and expressed in Lemma \ref{unit_interv} below, any number in $ \varOmega \left( n \right)$ i.e. any odd composite number  will appear within one integers interval, $ \varOmega_k.$ A trivial proof follows. 
\end{itemize}
Let us consider any odd composite number $p$: $ p \geq 9 \ \mbox{and odd means that } \  p = 3*3 + 2k, \ \mbox{for some} \ k \in \mathbb{N}.$ 
On the other hand, $  p  \leq p + 4 k \, \ k \in \mathbb{N} \ \mbox{given above means that} \ p \leq 3 \left( 3 + 2 k   \right).$ 
This means that $ p \in  \llbracket 9 \ .. \   3 \left( 3 + 2 k   \right) \rrbracket \  \mbox{for some given} \ k \in  \mathbb{N}: $ the latter is an interval of length $6k.$ Such interval has $k$ sub-intervals, each of length 6 with $p$ in one them. Each of these k sub-intervals of length 6 have boundaries $N_{k+1} = 3 \left( 3 + 2 k   \right) $   and $ N_k = N_{k+1} - 6  \Leftrightarrow N_k = 3 \left( 1 + 2 k   \right), \ k \in\mathbb{N^*} $ since the smallest odd composite number is $9$.
\begin{eqnarray}
\nonumber 
 N_k &=& N_{k+1} - 6 \nonumber \\
 & = & 3 \left( 1 + 2 k   \right) \nonumber \\
 & = & 3 \left( 3 + 2 l   \right) \ \mbox{for} \ l \in  \mathbb{N^*} \ \mbox{since } \ \left( 1 + 2 k   \right) \ \mbox{is odd. Then repeat for $N_{k - 1} = N_k - 6$.}
 \nonumber 
\end{eqnarray}
Since each boundary is an odd number multiple of 3, and $ \llbracket 9 \ .. \   3 \left( 3 + 2 k   \right) \rrbracket $ is finite, repeating the reasoning above establishes that $ p \in  \llbracket 3 \left( 1 + 2 m   \right)  \ .. \   3 \left( 3 + 2 m   \right) \rrbracket \ \mbox{for} \ m  \in  \mathbb{N^*}. $ We then conclude that any odd composite number $p$ is within an interval (of integers) $ \varOmega_k = \llbracket N_k \ .. \ N_{k+1}   \rrbracket $ of length 6 and whose boundaries are each a multiple of 3.

Let us then investigate what happens within each interval (of integers)  $ \varOmega_k:$ any pair of {\it consecutive} odd composite numbers has to be within one $ \varOmega_k,$, otherwise there is no way they could be consecutive (we refer the reader to the definition of consecutive odd composite number in Definition 2).
 \begin{lemm} \label{unit_interv}
  For any pair $ \left\lbrace p_1 \ ,  p_2 \right\rbrace $ of consecutive odd composite numbers each multiple of $3$, there always exist two and only two odd numbers in between.
 \end{lemm}
 A trivial proof consists in writing down any such pair with $p_1 = 3\left(2n + 1 \right)$ and $p_2 = 3\left(2n + 3 \right)$, $n \in \mathbb{N}.$ Since they are distant by 6, there are only two odd numbers in between, $p_{11}$ and $p_{12}$. Therefore only three possibilities for their separation: $ \left \vert   p_{11} - p_{12} \right \vert = $ $2$ or $4$ or $6$.  The trivial proofs for the following lemma consist in writing down, each time, the two consecutive odd composite numbers and identifying the odd numbers in between. 
 \begin{lemm} \label{one_int_only}
 For any pair $ \left\lbrace N_1, \ N_2 \right\rbrace $ of consecutive odd composite numbers, If $N_2 - N_1 = 2$ then there is one and only one integer between $N_1$ and $N_2$ and that integer is even. 
 \end{lemm}
 \begin{lemm} \label{one_pn_only}
 For any pair $ \left\lbrace N_1,\ \ N_2 \right\rbrace $ of consecutive odd composite numbers. If $N_2 - N_1 = 4$ then there exists one and only one prime number $p_1$ between $N_1$ and $N_2$. That prime number is given by $ p_1 = \frac{N_1 + N_2}{2} $. It is an Isolated prime.
 \end{lemm}
 \begin{lemm} \label{two_pn_only}
 For any pair $ \left\lbrace N_1, \ \ N_2 \right\rbrace $ of consecutive odd composite numbers. If $N_2 - N_1 = 6$ then there exists two and only two prime numbers $ p_1 $ and $ p_2 $ between $N_1$ and $N_2$. These prime numbers given by $ p_1 = \frac{N_1 + N_2}{2} - 1 $ and $ p_1 = \frac{N_1 + N_2}{2} + 1 $. $p_1$ and $p_2$ are Twin primes.
 \end{lemm}
\end{enumerate}
\label{proof_decision} 
\end{proof}

The eleven relations (tables), previously mentioned, yield the set of consecutive odd composite numbers below any integer $n \geq 9$. The three-criteria test in Theorem \ref{decision_test} applied to differences of consecutive pairs of the previous set yields $\mathcal{P}(n),$ the set of primes below $n,$ with the trivial ones $ \left\lbrace 2, 3, 5, 7 \right\rbrace $ added. In a nutshell, that is our novel deterministic scheme. 

\section{A Few Consequences}

In all subsequent sections, our preassigned threshold $m = 3 \left( 2 n + 1 \right)$ for some $ n \in  \mathbb{N}^*$. In other words, $m$ is a odd composite number multiple of $3$. This has the merit to lead to versatile relations exploiting the fact that the intervals (of integers) $ \Omega_k = \llbracket 3 \left( 1 + 2 k   \right)  \ .. \   3 \left( 3 + 2 k   \right) \rrbracket , k \in  \mathbb{N}^* $ provide a partitioning for $ \varOmega_{cocoon} \left( m \right) $, as illustrated in Figure \ref{fig:line_of_numbers}.

\subsection{Number of Primes Below Any $m = 3 \left( 2 n + 1 \right), n \in  \mathbb{N}$}

\begin{kor}
 Theorem \ref{decision_test} establishes, for the first time, that
 \begin{equation}
\pi(m)  = 4 + \left \vert A_4(m) \right \vert + 2 \left \vert  A_6(m) \right \vert  \ \mbox{for any $m = 3 \left( 2 n + 1 \right)$ where $ n \in  \mathbb{N}^*,$ } \\
\label{pi_of_m}
 \end{equation}
 \end{kor}
 
where  $ A_4(m) $ and $ A_6(m) $ are the sets of pairs of consecutive composites odd numbers distant by $4$ and $6$, respectively. The intuition of the relatively trivial proof for relation \ref{pi_of_m} is given by two facts. First, we trivially have $ \mathbb{N} \subseteq \llbracket 0 \ .. \ 8  \rrbracket \bigcup \displaystyle \bigcup^{\infty}_{k = 1, \  k \in \mathbb{N} } \varOmega_k, $ with $\varOmega_k $ as above. Adding $ \llbracket 0 \ .. \ 8  \rrbracket $ to it,  Figure \ref{fig:line_of_numbers} provides a robust visualization of the latter. Secondly, based on Theorem \ref{decision_test}, one realizes that $\p(m),$ the set of all Primes below $m$ is given by $\p(m) = \left \lbrace 2 , 3, 5 , 7 \right \rbrace  \bigcup $   $ \left \lbrace \mbox{ Primes within all pairs in $A_6(m)$ } \right \rbrace  \bigcup $ $ \left \lbrace \mbox{Primes within all pairs NOT in $A_6(m)$}  \right \rbrace $ where $  \left \lbrace \mbox{Primes within all pairs NOT in $A_6(m)$ } \right \rbrace =  \left \lbrace \mbox{Primes within all pairs in $A_2(m)$ } \right \rbrace  \bigcup  \left \lbrace \mbox{Primes within all pairs in $A_4(m)$ } \right \rbrace ,$ {\it for any $ m = 3 \left( 2n + 1 \right), \ n \in \mathbb{N^*}. $ All disjoint.} In other words, $ \lbrace 2 , 3, 5 , 7 \rbrace $, $ A_6(m) $, $ A_4(m) $ form a partition of $\p(m)$. With the latter, the Cardinal of $\p(m)$ (= $\pi(m)$) is straightforward by counting the Primes in each partition. Adding the latter yields relation \ref{pi_of_m}.
\begin{kor} Theorem \ref{decision_test} implies that for $m$ as above,
\begin{itemize}
\item $ A_6(m) $ corresponds to the set of Twin Primes below $m$ i.e to pairs of Primes distant by $2,$ below $m$;
\item $ A_4(m) $ is the set of primes below $m$ that each appear between a pair of COCOONS distant by $4$: we name them ''Isolated Primes";
\item $ A_2(m) $ is the set of pairs of COCOONS distant by $2$, below $m.$
\label{Polignac}
\end{itemize}
\end{kor}
Based on Euclid's Theorem, Equation \ref{pi_of_m} establishes that there are infinitely many "Isolated  primes" {\bf or} infinitely many Twin primes. The last two claims are major conjectures that we explore in the next sections.

\subsection{ $ \left \vert  A_2(m)  \right \vert  $ : Number of Pairs of Consecutive odd composite Numbers Below Any $m = 3 \left( 2 n + 1 \right), n \in  \mathbb{N^*}$, and distant by $2$}
For any $m = 3 \left( 2 n + 1 \right)$ with $ n \in  \mathbb{N}^*$,
we start by re-expressing, $\Omega_{even} \left( m \right),$ the set of all even numbers below $m$ as
\begin{eqnarray}
 \Omega_{even} \left( m \right) & = & \lbrace 0, 2, 4, 6, 8 \rbrace \bigcup^{M }_{k = 1, \  k \in \mathbb{N} } \lbrace \mbox{All even numbers in} \   \varOmega_k, \mbox{ with $ \varOmega_k = \llbracket N_1(k)  \ .. \   N_2(k) \rrbracket  $ } as \ above \rbrace \nonumber \\ 
& = & \lbrace 0, 2, 4, 6, 8 \rbrace \bigcup^{M }_{k = 1, \  N_1(k), \ N_2(k) \ COCOONS} \i_k  \bigcup^{M }_{j = 1, \  N_1(j), \ N_2(j) \ NOT \ Consecutive} \i_j  
\end{eqnarray}
where $ \i_k \equiv \lbrace \mbox{All even numbers in} \   \varOmega_k, \mbox{ with $ \varOmega_k = \llbracket N_1(k)  \ .. \   N_2(k) \rrbracket  $ } \rbrace .$

Based on the relation above, it is relatively trivial to establish that
$$ \Omega_{even} \left( m \right) = \lbrace 0, 2, 4, 6, 8 \rbrace \bigcup \lbrace  
\mbox{All evens in $   A_6(m)  $ } \rbrace \bigcup \lbrace  
\mbox{All evens in $   A_4(m)  $ } \rbrace \lbrace  
\mbox{All evens in $   A_2(m)  $ } \rbrace   . $$ 

Counting them yields 
\begin{eqnarray}
   \left \vert  \Omega_{even} \left( m \right) \right \vert  & = & 5 + 3\left \vert  A_6(m)  \right \vert +  2\left \vert  A_4(m)  \right \vert +  \left \vert  A_2(m)  \right \vert \\
   \frac{ m + 1}{2} & = & 5 + 3\left \vert  A_6(m)  \right \vert +  2\left \vert  A_4(m)  \right \vert +  \left \vert  A_2(m)  \right \vert . \nonumber
\end{eqnarray}
Figure \ref{fig:line_of_numbers} provides a robust visual help in the proof above. 
We have established that
\begin{equation}
    \left \vert  A_2(m)  \right \vert  =  \frac{m - 9}{2} - 3\left \vert  A_6(m)  \right \vert - 2\left \vert  A_4(m)  \right \vert.
   \label{A2m} 
\end{equation}
\subsection{Number of Consecutive Odd Composite Numbers Below Any $m \in  \mathbb{N^*}$: $\t(m)$}

For any $m = 3 \left( 2 n + 1 \right)$ with $ n \in  \mathbb{N}^*$, $\t(n),$ designates the number of consecutive odd composite numbers (COCOON) below any preassigned limit $ m:  \t(m) = \left \vert \varOmega_{COCOON} \left( m \right)  \right \vert .$ By construction  $ \varOmega_{COCOON} \left( m \right) $ is a finite set made of $ \left \lbrace u_i \right \rbrace_{i=1,2,3,...,l+1} .$ Let us introduce $ \tilde{\varOmega}_{COCOON} = \left \lbrace \left(u_1, u_2 \right),\left(u_2, u_3 \right), \left(u_3, u_4 \right), ..., \left(u_{l-1}, u_l \right), \left(u_l, u_{l + 1} \right) \right \rbrace .$  Clearly, $ \left \vert \tilde{\varOmega}_{COCOON}   \right \vert = l $ and  $ \left \vert \varOmega_{COCOON} \left( m \right)  \right \vert = 1 + l .$ In other words: $ \left \vert \varOmega_{COCOON} \left( n \right)  \right \vert = 1 + \left \vert \tilde{\varOmega}_{COCOON}   \right \vert .$ 
We now introduce $$ \tilde{\varOmega}^{(2)}_{COCOON} = \left \lbrace \left(u_i, u_j \right), \ \mbox{such that $u_j - u_i = 2$}   \right \rbrace, $$
$$ \tilde{\varOmega}^{(4)}_{COCOON} = \left \lbrace \left(u_i, u_j \right), \ \mbox{such that $u_j - u_i = 4$}   \right \rbrace, $$ 
$$ \tilde{\varOmega}^{(6)}_{COCOON} = \left \lbrace \left(u_i, u_j \right), \  \mbox{such that $u_j - u_i = 6$}   \right \rbrace.$$
Embedded in Theorem \ref{decision_test} is the fact that $$ \tilde{\varOmega}_{COCOON} = \tilde{\varOmega}^{(2)}_{COCOON} \bigcup \tilde{\varOmega}^{(4)}_{COCOON} \bigcup \tilde{\varOmega}^{(4)}_{COCOON} ,$$ with $$ \tilde{\varOmega}^{(2)}_{COCOON} \bigcap \tilde{\varOmega}^{(4)}_{COCOON} \bigcap \tilde{\varOmega}^{(4)}_{COCOON} = \o.$$
Therefore $$\left \vert \tilde{\varOmega}_{COCOON} \right \vert = \left \vert \tilde{\varOmega}^{(2)}_{COCOON} \right \vert + \left \vert \tilde{\varOmega}^{(4)}_{COCOON} \right \vert + \left \vert \tilde{\varOmega}^{(6)}_{COCOON} \right \vert,$$ hence  
\begin{thm}
 \begin{eqnarray}
 \nonumber 
  \left \vert \varOmega_{COCOON} \left( m \right)  \right \vert &=& 1 + \left \vert \tilde{\varOmega}_{COCOON} \right \vert , \\ \nonumber
  &=& 1 + \left \vert \tilde{\varOmega}^{(2)}_{COCOON} \right \vert + \left \vert \tilde{\varOmega}^{(4)}_{COCOON} \right \vert + \left \vert \tilde{\varOmega}^{(6)}_{COCOON} \right \vert, \\ \nonumber 
  \t(m)&=& 1 + \left \vert A_2(m) \right \vert + \left \vert A_4(m) \right \vert + \left \vert  A_6(m) \right \vert \ \forall \ m=3*\left(2n + 1  \right), \ n \ \ in \  \mathbb{N}^*.  
  \nonumber
 \end{eqnarray}
\label{no_of_oc}
\end{thm}
From Theorem \ref{decision_test}, we get that $\left \vert  A_2(n)  \right \vert$ is the number of even numbers below $n$. Each term in Theorem \ref{no_of_oc} is trivially computed via our deterministic scheme. 

\subsubsection{Some Intermediate Relations}

For any $ n \in \mathbb{N}, $ we now introduce the following sets:
\begin{enumerate}[(1)]
\item $ \varOmega_n = \left \lbrace 0, 1, 2, 3, ..., k, ..., n \ \mbox{with $n > 0$} \right \rbrace, $ the set of integers smaller than $n$;
\item $ \varOmega_{even} \left( n \right) = \left \lbrace 2k \leq n, \ k \in  \mathbb{N} \right \rbrace  $ the set of even integers smaller than $n$;  
\item $ \varOmega_{odd} \left( n \right) = \left \lbrace 2k + 1 \leq n, \ k \in  \mathbb{N} \right \rbrace  $ the set of odd integers smaller than $n$.
\item $ \varOmega_{c} \left( n \right)  = \left \lbrace \mbox{$p,$ an odd composite number with} \ p \leq n, \right \rbrace  $ the odd composite numbers smaller than $n$;  
\item $ \p \left( n \right) = \left \lbrace \mbox{$p,$ a prime number with} \ p \leq n \right \rbrace  $ the set of prime numbers smaller than $n$;
\end{enumerate}
We have
\begin{eqnarray}
\nonumber 
 \left \vert  \varOmega_n \right \vert &=&  \left \vert  \varOmega_{even} \left( n \right) \right \vert  + \left \vert  \varOmega_{odd} \left( n \right) \right \vert \\ \nonumber
 &=& n + 1.
\end{eqnarray}
If $n$ is odd then $$ \left \vert \varOmega_{even} \left( n \right) \right \vert = \frac{n + 1}{2} .$$
If $n$ is even then $$ \left \vert \varOmega_{even} \left( n \right) \right \vert = \frac{n}{2} + 1 .$$
Within $ \mathbb{N}$, 
\begin{eqnarray}
\nonumber
 \left \vert  \varOmega_{odd}\left( n \right) \right \vert &=& 1 + \left \vert  \varOmega_{coon} \left( n \right) \right \vert +  \left \vert  \varOmega_{primes} \left( n \right) \right \vert, \\ \nonumber
\label{all_odds}
 \end{eqnarray}
Therefore 
\begin{kor}
\begin{eqnarray}
\nonumber 
 \left \vert  \varOmega_n \right \vert &=&  \left \vert  \varOmega_{even} \left( n \right) \right \vert  +  \left \vert  \varOmega_{odd}\left( n \right) \right \vert     \\ \nonumber 
 n + 1 &=& \frac{n}{2}  + 1 + \left \vert  \varOmega_{coon} \left( n \right) \right \vert +  \left \vert  \varOmega_{primes} \left( n \right) \right \vert \ \ \mbox{If $n$ is even,} \\ \nonumber 
 n + 1 &=& \frac{n}{2}  + 1 + \t(n)  + \pi(n) \ \\ \nonumber 
 \end{eqnarray}
and
\begin{eqnarray}
\nonumber 
 \left \vert  \varOmega_n \right \vert &=&  \left \vert  \varOmega_{even} \left( n \right) \right \vert  +  \left \vert  \varOmega_{odd}\left( n \right) \right \vert     \\ \nonumber 
 n + 1 &=& \frac{n + 1 }{2} + \left \vert  \varOmega_{coon} \left( n \right) \right \vert +  \left \vert  \varOmega_{primes} \left( n \right) \right \vert  \ \ \mbox{If $n$ is odd.} \\ \nonumber 
 n + 1 &=& \frac{n + 1 }{2} + \t(n)  + \pi(n) \ \\ \nonumber 
 \end{eqnarray}
\label{up_to_n}
\end{kor}
In the above, one has to make sure not to double-count $2,$ the only integer both even and prime. 
Therefore,
\begin{kor}
 \begin{enumerate}[(1)] For any $ n \in \mathbb{N} $
  \item  $ \pi(n) = \frac{n}{2} - \t(n)$ If $n$ is even,
  \item  $ \pi(n) = \frac{n + 1}{2} - \t(n)$ If $n$ is odd.
 \end{enumerate}
 \label{count_them}
\end{kor}

\subsection{Revisiting The Twin Primes Conjecture}

Twin primes are pairs of primes which differ by two. It has been conjectured that there are infinitely many twin primes \cite{dePolignac1849}. Here we revisit this conjecture under a new light provided by our novel Deterministic scheme.  {\it Using the fact in Corollary \ref{Polignac} that $ A_6(n) $ corresponds to the set of Twin Primes below $m,$ our aim is to establish a lower bound for $ \left \vert A_6(n) \right \vert,$ of polynomial type in $m$ and conclude}. Here again, as in all the next sections, $m = 3 \left( 2 n + 1 \right)$ for some $ n \in  \mathbb{N}^*$: $m$ is an odd composite number multiple of $3.$ For such $m,$ odd integer, Corollary \ref{count_them} establishes that
\begin{equation}
\pi(m) + \t(m) = \frac{m + 1}{2}, \ \mbox{as expected}.
\label{eq3}
\end{equation}
Substituting $ \left \vert  A_2(m)  \right \vert  $ from Eqn \ref{A2m} into the expression for $ \t(m) $ in Theorem \ref{no_of_oc} yields
\begin{equation}
 4 \left \vert  A_6(m)  \right \vert + 7  = m -  2 \left(  \t(m) +  \left \vert  A_4(m)  \right \vert \right).
\label{twins}
\end{equation}
Equation \ref{twins} establishes that $ 0 < \frac{ 4\left \vert A_6(n) \right \vert + 7 }{m} = 1 - \frac{2}{m} \left(  \t(m) +  \left \vert  A_4(m)  \right \vert \right) , \ \mbox{for any $m$ as above}. $
Hence
\begin{eqnarray}
\label{sm}
 0 < 1 - \frac{2}{m}  \left(  \t(m) +  \left \vert  A_4(m)  \right \vert \right) < 1 \ \mbox{as $ ( \left(  \t(m) +  \left \vert  A_4(m)  \right \vert \right)  > 0),$ \ for any $m$ as above.}
\end{eqnarray}
Let us introduce then $$ \s = \left\lbrace 1 - \frac{2}{l} \left(  \t(l) +  \left \vert  A_4(l)  \right \vert \right) \ \mbox{such that} \ l = 3 \left( 2 k + 1 \right), \ k \in  \mathbb{N}^* \right\rbrace .$$
Due to the relation in \ref{sm}, the non empty set $\s$ is bounded, hence $Inf \left( \s \right)$ exists, is unique and finite. Since $\s$ does not depend on $m$ and $Inf \left( \s \right)$ is unique, the latter does not depend on $m$. We observe that for any $\alpha > 0$ such that 
\begin{equation}
 0 < \alpha  \leq 1 -  \frac{2 }{m} \left(  \t(m) +  \left \vert  A_4(m)  \right \vert \right) , \ \mbox{for any m as above,}
\label{alpha}
\end{equation}
one gets
\begin{equation}
 m - 2 \left(  \t(m) +  \left \vert  A_4(m)  \right \vert \right)  > \alpha m \Leftrightarrow m - 2 \left(  \t(m) +  \left \vert  A_4(m)  \right \vert \right) = 4 \left \vert  A_6(m)  \right \vert + 7  > \alpha m, \ \mbox{with $ \alpha > 0 $ and m as above.}
\label{ineq} 
\end{equation}
{\it Finding a suitable $\alpha > 0$ satisfying the condition in the relation \ref{alpha} would therefore provide us with a lower bound of polynomial type, in $m$, for $ \left \vert A_6(m \right \vert $, as in \ref{ineq}. Recalling that $ \left \vert A_6(m) \right \vert $ is the number of Twin Primes below $m,$ would then help to conclude when
$ m \longrightarrow + \infty.$} Let us therefore find such real number $\alpha > 0.$

In general, $ Inf \left( \s \right) \geq 0.$ However, due to the relation \ref{sm} and by its definition, the only option is $ Inf \left( \s \right) > 0.$ The latter is also guaranteed by the fact that $\mathbb{Q}$ is dense in $\mathbb{R}.$ We then introduce
\begin{equation} 
\alpha = Inf \left( \s \right). 
\end{equation}
{\it $\alpha$ is independent of $m$ } and we have $$ 0 < \alpha \leq 1 - \frac{ 2 }{m} \left(  \t(m) +  \left \vert  A_4(m)  \right \vert \right)  , \ \mbox{for any $m$ as above.}$$
For the chosen $ \alpha > 0,$
$$ 0 < \alpha \leq 1 -  \frac{ 2 }{m} \left(  \t(m) +  \left \vert  A_4(m)  \right \vert \right)  \Leftrightarrow  m - 2 \left(  \t(m) +  \left \vert  A_4(m)  \right \vert \right) = 4 \left \vert A_6(n) \right \vert + 7  > \alpha m, $$
for any $m = 3 \left( 2 n + 1 \right)$ \ with $ \ n \in  \mathbb{N}^* \ , and \ \alpha > 0 .$  In other words, we have established that 
$$  2\left \vert A_6(n) \right \vert > \frac{\alpha m - 7} {2}, \ \mbox{for any} \ m = 3 \left( 2 n + 1 \right) \ \mbox{with} \  n \in  \mathbb{N}^*, \ \mbox{and} \ \alpha > 0 .$$
Since  $$ \displaystyle \lim_{ m \longrightarrow \infty } \left( \frac{\alpha m - 7}{2} \right) = + \infty, $$ we deduce that
$$  \displaystyle \lim_{ m \longrightarrow \infty } 2 \left \vert A_6(m) \right \vert  = + \infty, $$
where we know from a previous section that $ 2 \left \vert  A_6(m)  \right \vert $ corresponds to the cardinal of the set of Twin Primes.
Through $ \pi(m), $ in Equation \ref{pi_of_m}, the latter provides a new proof for Euclid's Theorem. Through $\t(m),$ Theorem \ref{no_of_oc} establishes that there is an infinite number of odd composite numbers. {\it More importantly, the last relation establishes, for the first time, that there are infinitely many Twin Primes.}

\subsection{There Are Infinitely Many Isolated Primes}
We recall that $ A_4(m) $ is the set of primes below $m$ that each appear between a pair of COCOONS distant by $4$: they are named Isolated Primes. We established in a previous section that
\begin{equation}
    \left \vert  A_2(m)  \right \vert  =  \frac{m - 9}{2} - 3 \left( \left \vert  A_6(m)  \right \vert - 2 \left \vert  A_4(m)  \right \vert \right).
   \label{AA2m} 
\end{equation}

Substituting $ \left \vert  A_2(m)  \right \vert  $ from Eqn \ref{AA2m} into the expression for $ \t(m) $ in Theorem \ref{no_of_oc} yields
\begin{equation}
 2 \left \vert  A_4(m)  \right \vert + 7  = m -  2 \left(  \t(m) + 2 \left \vert  A_6(m)  \right \vert \right).
\end{equation}
In other words, 
\begin{equation}
0 < \frac{2 \left \vert  A_4(m)  \right \vert + 7}{m}  = 1 -  \frac{2}{m} \left(  \t(m) + 2 \left \vert  A_6(m)  \right \vert \right) < 1.
\label{isolated}
\end{equation}
As in the previous section, we introduce $$ \tilde{\s} = \left\lbrace 1 - \frac{2}{l} \left(  \t(l) + 2 \left \vert  A_6(l)  \right \vert \right) \ \mbox{such that} \ l = 3 \left( 2 k + 1 \right), \ k \in  \mathbb{N}^* \right\rbrace .$$
$ \tilde{\s} $ does not depend on $m.$ Due to the relation in \ref{isolated}, the non empty set $\tilde{\s}$ is bounded, hence $Inf \left( \tilde{\s} \right)$ exists, is unique and finite. $Inf \left( \tilde{\s} \right)$ does not depend on $m,$ since $\tilde{\s}$ does not depend on $m$ and $Inf \left( \tilde{\s} \right)$ is unique. 

In general, $ Inf \left( \tilde{\s} \right) \geq 0.$ However, due to the relation \ref{sm} and by its definition, the only option is $ Inf \left( \tilde{\s} \right) > 0.$ The latter is also guaranteed by the fact that $\mathbb{Q}$ is dense in $\mathbb{R}.$ We then introduce
\begin{equation} 
\tilde{\alpha} = Inf \left( \tilde{\s} \right) \ \mbox{is independent of m.}  
\end{equation}
We then have $$ 0 < \tilde{\alpha} \leq 1 - \frac{ 2 }{m} \left(  \t(m) + 2  \left \vert  A_6(m)  \right \vert \right)  , \ \mbox{for any $m$ as above.}$$
For the chosen $ \tilde{\alpha} > 0,$
$$ 0 < \tilde{\alpha} \leq 1 -  \frac{ 2 }{m} \left(  \t(m) + 2  \left \vert  A_6(m)  \right \vert \right)  \Leftrightarrow  m - 2 \left(  \t(m) + 2 \left \vert  A_6(m)  \right \vert \right) = 4 \left \vert A_4(n) \right \vert + 9  \geq \tilde{\alpha} m, $$
for any $m = 3 \left( 2 n + 1 \right)$ \ with $ \ n \in  \mathbb{N}^* \ , and \ \tilde{\alpha} > 0 .$  In other words, we have established that 
$$  \left \vert A_4(n) \right \vert \geq \frac{ \tilde{\alpha} m - 9} {4}, \ \mbox{for any} \ m = 3 \left( 2 n + 1 \right) \ \mbox{with} \  n \in  \mathbb{N}^*, \ \mbox{and} \ \tilde{\alpha} > 0 .$$
Since  $$ \displaystyle \lim_{ m \longrightarrow \infty } \left( \frac{ \tilde{ \alpha} m - 9}{4} \right) = + \infty, $$ we deduce that
$$  \displaystyle \lim_{ m \longrightarrow \infty } \left \vert A_4(m) \right \vert  = + \infty, $$
where we know from a previous section that $ \left \vert  A_4(m)  \right \vert $ corresponds to the cardinal of the set of "Isolated" Primes, below m.
The last relation above unambiguously proves, for the first time, that {\it there are infinitely many "Isolated" Primes}.

\section{\textbf{Conclusion \& Future Works }}

A new deterministic scheme for characterizing the organization of primes has been established. It yields the set of consecutive odd composite numbers below any preassigned limit $n \in\mathbb{N}.$ $\p(n),$ the set of all primes below $n$ is trivially deduced from the latter. At the core of the scheme are eleven generic tables, coupled with a three-criteria test applied on differences between pairs of consecutive odd composite numbers obtained from the tables. $\left \vert \p(n) \right \vert \equiv  \pi(n).$ We establish that  $\pi(m) = 4 + \left \vert A_4(m) \right \vert + 2 \left \vert  A_6(m) \right \vert  \ for \ m = 3 \left( 3l + 1  \right), \ \forall \  l \in \  \mathbb{N}^*.$ 
$ \left \vert  A_4(m)  \right \vert $ and $ \left \vert  A_6(m)  \right \vert$ are the cardinals of the {\it sets of pairs of consecutive odd composite numbers}, below $m$ and having, each, a difference of $4$ and $6,$ respectively. $  \left \vert  A_6(n) \right \vert$ and $ \left \vert A_4(n) \right \vert $ correspond to the number of Twin primes and `` Isolated'' primes, below any $m$ respectively. $\t(m)$ is the number of consecutive odd composite numbers below $m$.We establish that $\t(m) = 1 + \left \vert A_2(m) \right \vert + \left \vert A_4(m) \right \vert + \left \vert  A_6(m) \right \vert$ where $ \left \vert  A_2(m)  \right \vert $ is the number of pairs of consecutive odd composite numbers, below $m$, and distant by $2$. {\it We prove that $ A_6(m) $ is the set of Twin Primes below $m,$ $A_4(m) $ is the set of primes below $m$ that we call "isolated" as they each appear between a pair of consecutive odd composite numbers  distant by $4$.} $\t(m)$ is the number of consecutive odd composite numbers below $m$. We establish that $\t(m) = 1 + \left \vert A_2(m) \right \vert + \left \vert A_4(m) \right \vert + \left \vert  A_6(m) \right \vert$ where $ \left \vert  A_2(m)  \right \vert $ is the number of pairs of consecutive odd composite numbers, below $m$, distant by $2$. We prove that {\it $ \left \vert  A_2(m)  \right \vert  = \frac{m - 9}{2} - 3 \left \vert  A_6(m)  \right \vert - 2 \left \vert  A_4(m)  \right \vert .$ } With {\it $m$ an odd number}, combining $\pi(m)$, $ \t(m),$ and $ \left \vert  A_2(m)  \right \vert, $ we establish that  $$ 4 \left \vert  A_6(m)  \right \vert + 7  = m -  2 \left(  \t(m) +  \left \vert  A_4(m)  \right \vert \right)  > 0, \ hence \  0 < 1 -  \frac{ 2 }{m} \left(  \t(m) +  \left \vert  A_4(m)  \right \vert \right)  < 1 , \ \mbox{for any $m$ as above.} $$ The non empty set $\s \equiv \left\lbrace 1 -  \frac{ 2 }{l} \left(  \t(l) +  \left \vert  A_4(l)  \right \vert \right) \ \mbox{such that} \ \mbox{l = 3(2k + 1), $k \in \mathbb{N}^*$} \right\rbrace \subset \mathbb{R} $ is bounded and independent of m. 
Hence $$ Inf \left( \s \right)  \ \mbox{{\it exists, is unique, finite and independent of m.}}$$ 
By definition, $ Inf \left( \s \right) > 0.$ $\mathbb{Q}$ Dense in $\mathbb{R}$ also guarantees the latter. We then introduce
\begin{equation} 
\alpha = Inf \left( \s \right): \ \mbox{independent of m.} 
\end{equation}
We have $$ 0 < \alpha \leq 1 -  \frac{ 2 }{m} \left(  \t(m) +  \left \vert  A_4(m)  \right \vert \right) , \ \mbox{for any m = 3(2n + 1), $ n \geq 1$, an integer.}$$
In others words, $$ m - 2 \left( \t(m) +  \left \vert  A_4(m)  \right \vert \right) = 4 \left \vert  A_6(m)  \right \vert + 7 \geq \alpha m \Leftrightarrow 2 \left \vert  A_6(m)  \right \vert \geq \frac{\alpha m - 7}{2} , \ \mbox{ with $ \alpha > 0$ and m as above.} $$
{\it The last relation establishes, for the first time, that $$ \displaystyle \lim_{ m \longrightarrow \infty }  2 \left \vert  A_6(m)  \right \vert  = + \infty , \ \mbox{with  $ 2 \left \vert  A_6(m)  \right \vert, $ the cardinal of the set of Twin Primes below $m.$ } $$ }
Adopting a reasoning similar to the one above, we have established that $$ \displaystyle \lim_{ m \longrightarrow \infty }  \left \vert  A_4(m)  \right \vert  = + \infty \ \mbox{: there are infinitely many "Isolated" Primes}.$$
This also serves as a new proof of Euclid's Theorem. Through $\t(m),$ Equation in Theorem \ref{no_of_oc} establishes that there is an infinite number of odd composite numbers. More importantly the first of the two relations above unambiguously proves, for the first time, that {\it there is an infinite number of Twin Primes} and the second relation establishes for the first time too that {\it there are infinitely many "Isolated" Primes}. Coupling both results, it is established, {\it for the first time}, that {\it there are infinitely many pairs of primes distant by $2k$ when $k \in \left \lbrace 0,1 \right \rbrace$.} Significant results, collateral to the ones above will be the object of follow up publications. Among the latter are the establishment of novel entire functions of relative interest for the characterization of the organization of Primes. 
\newpage
{\bf Acknowledgements:}\\
PO acknowledges funding from the South African National Research Foundation (NRF). PO acknowledges the support and the generous hospitality provided by the ``L'Universit\'e Mahatma Gandhi'' in Conakry and the `` Minist\`ere de l'Enseignement Superieur, de la Recherche Scientifique et de l'Innovation '' of the Republic of Guinea during the research visits towards the present work. G.H acknowledges support from ``L'Universit\'e Mahatma Gandhi'' in Conakry and the `` Minist\`ere de la Recherche Scientifique, de la Recherche Scientifique et de l'Innovation '' of the Republic of Guinea, as well as from the State of Guinea. The authors wish to express their gratitude to Fode Idrissa Soumare for the critical support provided on computational aspects of the present work. Our words of thanks are also adressed to Bob Osano, at UCT, for hosting us.

\end{document}